\documentclass{amsart} \usepackage{amssymb,amsfonts,amsxtra}
\usepackage[all]{xy} \usepackage{fullpage}
\newtheorem{theorem}{Theorem}[section]
\newtheorem{cor}[theorem]{Corollary}
\newtheorem{lem}[theorem]{Lemma}
\newtheorem{prop}[theorem]{Proposition}

\newtheorem{example}[theorem]{Example}
\DeclareMathOperator{\Map}{Map} \DeclareMathOperator{\Hom}{Hom}
\DeclareMathOperator{\im}{im} \DeclareMathOperator{\ev}{ev}
\def\Num{\underline{Num}} 
 \def\lra{\longrightarrow}
  \def\ev{ev}
 \def\f{\mathbb{F}} 

\newtheorem{defi}[theorem]{Definition}
\newtheorem{rem}[theorem]{Remark} \numberwithin{equation}{section}
\thanks{This research was partially supported by the EPSRC grant
No. GR/R84276/01} \begin{document} \title[Hochschild cohomology
and moduli spaces...] {The discrete Gelfand transform and its
dual} \author{V.Buchstaber\&A.Lazarev} \address{Steklov
Mathematical Institute, Russian Academy of Sciences, Gubkina 8
Moscow 117966, Russia.}\email{buchstab@mendeleevo.ru
}\address{Mathematics Department, University of Bristol, Bristol,
BS8 1TW, England.} \email{a.lazarev@bristol.ac.uk}
\begin{abstract} We consider the transformation $\ev$ which
associates to any element in a
     K-algebra A a function on the the set of its K-points. This  is the
     analogue of the  fundamental Gelfand transform.
      Both $\ev$ and its dual $\ev^*$ are the maps from a discrete K-module to a
     topological K-module and we investigate in which case the image of each
     map is dense. The answer is nontrivial for various choices of K and A
     already for A=K[x], the polynomial ring in one variable.
     Applications to the structure of algebras of cohomology operations are given.
\end{abstract} \maketitle \section{introduction} Let $K$ be a
commutative ring and $A$ be a commutative $K$-algebra. The set of
$K$-points $A(K)$ is the set of $K$-algebra homomorphisms $A\lra
K$. The evaluation map is the canonical map associating to any
element $a\in A$ the function $\ev(a):A(K)\lra K$. For a $K$-point
$x:A\lra K$ the value of $\ev(a)$ at $x$ is $x(a)$. Thus we get a
map $\ev:A\lra \Map(A(K),K)$. Its $K$-dual is a map
$\ev^*:K[A(K)]\lra A^*$ where $K[A(K)]$ is the free $K$-module on
the set $A(K)$. The study of the evaluation map and its dual
(defined in a slightly more general way) is the main purpose of
this paper. Note that the evaluation map is a direct analogue of
the well-known Gelfand transform which identifies a commutative
$C^*$-algebra with the algebra of functions on the set of its
maximal ideals. We do not know if the analogue of the dual
evaluation map has been considered in the context of Banach
algebras.

Both $\ev$ and $\ev^*$ are maps from a discrete $K$-module to a
topological $K$-module and we investigate in which case the image
of each map is dense. The answer is nontrivial already for
$A=K[x]$, the polynomial ring in one variable.

Informally speaking, there are obstructions to $\ev$ and $\ev^*$
having a dense image related to the existence of \emph{numerical
functions} and \emph{numerical functionals}. In the simplest case
$K=\mathbb{Z}$, $A=K[x]$ numerical functions are known by the name
of \emph{numerical polynomials}. A polynomial in one variable with
rational coefficients is called numerical if it takes integer
values evaluated at integers. Numerical polynomials appeared in
algebro-topological
literature \cite{Ba}, \cite{Hub} and
\cite{Eke}. There is also a large body of literature in which numerical polynomials
(called usually \emph{integer-valued polynomials}) are studied from the purely algebraic
point of view. The standard reference is the monograph \cite{CC}.

The dual notion is that of a numerical functional. For
$K=L=\mathbb{Z}$ a numerical functional is a linear combination of
$\mathbb{Z}$-points of $A$ \emph{with rational coefficients} which
determines an integer functional on $A$. As far as we know the
concept of a numerical functional is new.

Theorem \ref{main1} gives a rather complete answer when the map
$\ev$ has a dense image. The corresponding question of the map
$\ev^*$ is more subtle and a partial answer is provided by Theorem
\ref{main2}.

We also give some applications of our results to the structure of
algebras of cohomology operations in mod $p$ homology and complex
cobordisms. Recall a well-known theorem of Morava, cf. \cite{Ra},
Theorem 6.2.3 which states that the Morava stabilizer algebra
$S(n)$ is isomorphic to the dual of a certain group ring after a
suitable extension of scalars. The algebra $S(n)$ is essentially
the even part of the cooperation algebra for $K(n)$, the $n$th
Morava $K$-theory. This theorem reduces the cohomology of $S(n)$
to that of the continuous cohomology of the corresponding group,
the so-called Morava stabilizer group. This is of great importance
because of the connection with the stable homotopy groups of
spheres via certain generalized Adams-Novikov spectral sequences.

It is natural to ask if there are other contexts in which algebras
of cohomology operations are related with group rings. In the
context of usual cohomology theory mod $p$ the analogue of $S(n)$
is $\mathcal{A}$, the full Steenrod algebra if $p=2$ and the
algebra of reduced powers if $p>2$. In the context of complex
cobordisms the relevant analogue is the so-called
Landweber-Novikov algebra $S$.

The above mentioned theorem of Morava turns on the fact that the
Morava stabilizer algebra
$S(n)=\f_{p^n}[t_1,t_2,\ldots,]/(t_i^{p^n}-t_i;i>0)$ is
\emph{semisimple} i.e. isomorphic to the infinite sum of copies of
$\f_{p^n}$. This is not true for either $S^*$ or $\mathcal{A}^*$
and therefore we cannot expect these algebras to be isomorphic to
the duals of some group rings.

Instead, our result states that after an appropriate extension of
scalars and completion the algebras $\mathcal{A}$ and $S$ contain
dense Hopf subalgebras isomorphic to group rings of certain groups
which we explicitly identify.

Note that in both cases we don't have the exact equality between
the completed algebras of operations and the corresponding group
rings. Because of that the relationship between their cohomologies
is far more subtle than in the Morava case. We intend to
investigate this relationship in a future work.

The present paper contains a detailed exposition of results
announced in \cite{BL}.
\\

{\bf Notation}. Throughout the paper $L$ will denote a fixed
commutative ring with a fixed infinite subring $K$ so that $L$
becomes a $K$-algebra. Furthermore, $L$ is supposed to be an
integral domain unless specifically said otherwise. The field of
fractions of $L$ is denoted by $\hat{L}$. The unadorned tensor
product $\otimes$ will stand for $\otimes_K$. For an $L$-module
$M$ the notation $\Hom_L(M,L)$ will be shortened to $M^*$. The set
of \emph{$K$-points} of a commutative $K$-algebra $A$ is denoted
by $A(K)$. In other words, $A(K)$ is the set of $K$-algebra maps
$A\lra K$. For two sets $S$ and $T$ the collection of all maps
from $S$ to $T$ is denoted by $\Map(S,T)$. The notation $L[S]$
stands for a free $L$-module on a set $S$ and a basis vector in
$L[S]$ corresponding to $s\in S$ is denoted by $[s]$. The prime
field consisting of $p$ elements is denoted by $\f_p$.

\section{The evaluation map and its dual: generalities}
Let $V$ be a free $L$-module. We say that a submodule $U\subset V$
is \emph{admissible} if $U$ is a free $L$-module of finite rank
and $U$ splits off $V$ as a direct summand in the category of
$L$-modules. Set $V^*:=\Hom_L(V,L)$. Then $V^*$ has the structure
of a topological $L$-module where the fundamental system of
neighborhoods of $0$ is formed by $L$-modules $V_U$:
\[V_U=\{l\subset V^*:l(v)=0 \mbox{~for all~} v\in U\}.\] Here
$U$ run through all free $L$-submodules of $V$ of finite rank.
Alternatively, we may require all submodules $U$ to be admissible.
Clearly the collection of admissible submodules is cofinal in the
collection of all finite rank free submodules, so we get the same
topology.

Now let $ V, W$ be two $L$-modules where $V$ is free and $W$ is
endowed with an $L$-linear topology.
\begin{defi} A $L$-module homomorphism $\phi:V\lra W$ is called an
\emph{asymptotically split} monomorphism if for any admissible
submodule $U\subset V$ the composition $U\hookrightarrow V\lra W$
admits a continuous left inverse $L$-module homomorphism $W\lra
U$.\end{defi}
\begin{rem}Note that an asymptotically split monomorphism is clearly a
monomorphism. In addition if $V$ is a free $L$-module of finite
rank then an asymptotically split monomorphism is simply a
(continuously) split monomorphism.\end{rem} Now let $V,W$ be free
$L$-modules and $f:W\lra V^*$ be a $L$-homomorphism. Its $L$-dual
homomorphism $f^*:V\lra W^*$  is defined as usual by the formula
$(f^*(v))(w)=(f(w))(v)$ where $v\in V$ and $w\in W$.
\begin{lem}\label{dense}Let $V,W$ be free $L$-modules. Then a $L$-homomorphism
$f:W\lra V^*$ has a dense image if and only if the dual
homomorphism $f^*:V\lra W^*$ is an asymptotically split
monomorphism.\end{lem}
\begin{proof}Suppose that $f$ has a dense image and let $U$ be any
admissible submodule in $V$. Then the definition of the linear
topology in $V^*$ directly implies that the composite map
\[W\lra V^*\lra V^*/V_U\] is an epimorphism. Since $U$ is
admissible the map $V^*\lra U^*$ dual to the inclusion $U\subset
V$ is onto and therefore $V^*/V_U\cong U^*$. Since $U^*$ is a free
$L$-module, the epimorphism $W\lra V^*/V_U\cong U^*$ is split by a
map $U^*\lra W$. The dual to the latter map is a continuous map
$W^*\lra U$ which is a left inverse to the composition
$U\hookrightarrow V\lra W^*$.

Conversely, suppose that $f^*$ be an asymptotically split
monomorphism and $U$ is an admissible submodule in $V$. The
$L$-dual to the composite map \[U\hookrightarrow V\lra W^*\] will
be the split epimorphism $W\lra U^*\cong V^*/V_U$ from which it
follows that $f$ has a dense image.\end{proof}
\begin{rem}If $L$ is a field then any monomorphism $V\lra W^*$ is asymptotically split.
Indeed, let $U\subset V$ be an admissible submodule in $V$. Since
the composition $U\hookrightarrow V\lra W^*$ is a monomorphism
there exists a subspace $W_{U^\prime}\subset W^*$ for which the
composition \[U\hookrightarrow V\lra W^*\lra W^*/W_{U^\prime}\] is
also a monomorphism. This last monomorphism necessarily splits and
determines a \emph{continuous} splitting of the inclusion
$U\hookrightarrow W^*$.
\end{rem}

Recall that $K$ is a fixed  infinite subring in $L$ and let $A$ be
a $K$-algebra which is free as a $K$-module. Define the evaluation
homomorphism $\ev_{L/K}$ as follows:
\begin{equation}\label{evaluation}\ev_{L/K}:A\otimes L\lra
\Map(A(K),L):\ev(a\otimes l)(x)=lx(a).\end{equation} Here $x:A\lra
K$ is an element in $A(K)$ and $l\in L$. Clearly, $\ev_{L/K}$ is a
homomorphism of $L$-modules.

The $L$-algebra $\Map(A(K),L)$ can be supplied with an $L$-linear
topology where the fundamental system of neighborhoods of $0$ is
formed by $L$-modules\[A_X=\{g\in \Map(A(K),L):g(x)=0 \mbox{~for
all~}x\in X\}.\] Here $X$ runs through all finite subsets in
$A(K)$.

Let us denote by $L[A(K)]$ the set of all continuous $L$-linear
homomorphisms from $\Map(A(K),L)$ into $L$. It follows that
$L[A(K)]$ is a free $L$-module generated by the set $A(K)$.
Moreover, the canonical homomorphism
\[\gamma:L[A(K)]^*\lra \Map(A(K),L):\gamma(s)(x)=s(x); ~~s\in L[A(K)]^*\]
is a homeomorphism.

Let us consider the homomorphism
\begin{equation}\label{coevaluation}\ev^*_{L/K}:L[A(K)]\lra (A\otimes L)^*\cong
\Hom_{K}(A,L)\end{equation} dual to $\ev_{L/K}$. The study of the
maps (\ref{evaluation}) and (\ref{coevaluation}) is our main goal.
\begin{rem}In the case when $A$ is a Hopf algebra over $K$ with
comultiplication $\Delta:A\lra A\otimes A$ the set $A(K)$ becomes
a group under the convolution product:
\[[x][y](a)=(x\otimes y)(\Delta(a));~ x,y\in A(K).\] Therefore
$L[A(K)]$ is in this case the group algebra of $A(K)$ over $L$ and
also a Hopf algebra over $L$. Furthermore, $A^*$ is a
\emph{topological} Hopf algebra. That means that $A^*$ is a
topological algebra and the diagonal in $A^*$ is a continuous map
$A^*\lra A^*\hat{\otimes}A^*:=(A\otimes A)^*$ from $A^*$ into the
completed tensor square of $A^*$. Similarly $\Map(A(K),L)$ is a
topological Hopf algebra over $L$ dual to $L[A(K)]$. Observe that
the $A$ and $L[A(K)]$ could also be considered as topological Hopf
algebras with the discrete topology.
\end{rem}
The following well-known fact is easily checked by direct
inspection of definitions.
\begin{lem} If $A$ is a Hopf algebra over $K$ then $\ev_{L/K}$ and $\ev^*_{L/K}$ are maps
of topological Hopf algebras over $L$.\end{lem}

The $L$-module $\Hom_{K}(A,L)$ has a structure of an $A$-module.
Indeed, let $a\in A$ and $t\in \Hom_{K}(A,L)$. Then $at\in
\Hom_{K}(A,L)$ is defined by $(at)b=t(ab)$ for $b\in A$.

Likewise the $L$-module $L[A(K)]$ admits a structure of an
$A$-module via the map
\[A\otimes L[A(K)]\lra K[A(K)]:a\otimes [x]=x(a)[x].\]
Here $x$ is a $K$-algebra map $A\lra K$ and $[x]$ is the
corresponding basis element in $L[A(K)]$.
\begin{rem} Of course, $A^*$ and $L[A(K)]$ are naturally \emph{right}
$A$-modules, but since $A$ is commutative there is no difference
between left and right $A$-modules.\end{rem}
\begin{example}
Let $A=K[z]$, the ring of polynomials on one generator $z$. Then
$\Hom_K(A,L)=\Gamma_L$ is the so-called \emph{divided power ring}
with coefficients in $L$. It can be identified with a topological
vector space over $L$ with topological basis $u_n$, where
\[u_n(z^k)=\delta^n_k.\] Here $\delta^n_k$ is the Kronecker symbol.
In this basis the action of $A$ on $\Gamma_L$ is specified by the
formula $z^k\cdot u_n=u_{n-k}$ where we assume that $u_i=0$ for
$i<0$. Clearly the set $A(K)$ is identified with $K$ and this
gives an isomorphism of $L$-modules $L[A(K)]\lra L[K]$ where
$L[K]$ is the free $L$-module generated by the elements of $K$.

The $L$-algebra map $\ev_{L/K}:L[z]\lra \Map(K,L)$ is given by the
formula $\ev_{L/K}(z)(k)=k$. The dual map $\ev^*_{L/K}:L[K]\lra
\Gamma_L$ is given by the formula \[ev^*_{L/K}[k]=\sum_{i\geq
0}k^iu_i.\]
\end{example}

We now have the following lemma whose proof is a direct inspection
of definitions:
\begin{lem}The map $\ev^*_{L/K}:L[A(K)]\lra \Hom_K(A,L)$ is a map of
$A$-modules.\end{lem} Next we have the following proposition whose
proof is similar to Artin's proof on linear independence of
characters, cf. \cite{lang}.
\begin{prop}\label{artin}The homomorphism $\ev^*_{L/K}$ is a monomorphism.\end{prop}
\begin{proof}Let $M= \ker \ev^*_{L/K}.$ It is an $L$-submodule in $L[A(K)]$.
According to the previous lemma $M$ is in fact an $A$-submodule.
Assume that $M \neq 0$.

As before we denote by $[x]$ the element in $L[A(K)]$
corresponding to a $K$-algebra map $x:A\lra K$. We will define the
map $\xi:M\lra \mathbb{Z}$ as follows. Write an arbitrary element
$s\in M\subset L[A(K)]$ as $s=\sum_{i=1}^N\alpha_i[x_i]$ where all
$x_i$ are pairwise distinct and all $\alpha_i$ are nonzero. Then
$\xi(s)=N$. Further set $m=\min_{s\in M}\xi(s)$. Clearly, $m>1$.
Consider an element $s\in M$ of the form
\[s=\alpha_1[x_1]+\alpha_2[x_2]+\ldots+\alpha_m[x_m].\]
Choose an element $a\in A$ for which $x_1(a)\neq x_2(a)$. Then
\[\hat{s}:=x_1(a)s-as=\sum_{i=2}^m\alpha_i(x_1(a)-x_i(a))[x_i]\]
is a nonzero element in $M$ such that $\xi(\hat{s})=m-1$. (Note
that here we used the fact that elements in $K$ are not zero
divisors in $L$.) This is a contradiction and our proposition is
proved.
\end{proof}
\begin{rem}
Without the condition that $K$ have no zero divisors in $L$ the
above result is no longer true. Indeed, let $A=K[x]$. Then for any
$k\in K$ we have
\[\ev^*_{L/K}([k]-0)=k(\sum_{i=1}^\infty k^{i-1}u_i)\in A^*= \Gamma_L.\]
Suppose that there exist elements $k\in K$ and $l\in L$ for which
$kl=0\in L$. It follows that the element $l([k]-[0])$ belongs to
the kernel of $\ev^*_{L/K}$.
\end{rem} Combining the above result with Lemma
\ref{dense} we get the following result:
\begin{cor}If $L$ is a field then the homomorphism $\ev_{L/K}:A\lra
\Map(A(K),L)$ has a dense image.\end{cor}
\begin{rem}If, in addition, the field $K$ is algebraically closed
and the $K$-algebra $A$ is finitely generated and has no nilpotent
elements then Hilbert's Nullstelensatz says that the $\ev_{L/K}$
is a monomorphism. In this case we conclude that $\ev^*_{L/K}$ is
a monomorphism with a dense image.

Later on we will consider the case $A=L[x]$ in some detail. It is
well-known and easy to see that in this case the map
$\ev_{L/K}:L[x]\lra \Map(L,K)$ is monomorphic. It also clear that
the condition that $K$ have no zero divisors cannot be
avoided.\end{rem}

We finish this section by considering the evaluation map and its
dual for tensor products of algebras. Let $A_i$ be $K$-algebras
where $i$ runs through an indexing set $I$, not necessarily
finite. Denote by $\ev_{L/K}(i)$ the evaluation map $A_i\otimes
L\lra \Map(A_i(K),L)$. Let $A$ be the tensor product of $A_i$;
more precisely \[A=\bigotimes_{i\in
I}A_i=\lim_{\rightarrow}\bigotimes_{k=1}^N A_{i_k}\] where the
limit is taken over all finite collections $(i_1,\ldots i_N)$ in
$I$. We will keep the notation $\ev_{L/K}$ for the evaluation map
$A\otimes L\lra \Map(A(K),L)$.

Similarly we denote the dual evaluation map $L[A_i(K)]\lra
\Hom_K(A_i,L)$ by $\ev^*_{L/K}(i)$. The dual evaluation map for
$A$ will still be denoted by $\ev^*_{L/K}$. Then we have the
following result.
\begin{prop}\label{tensor}
The map $\ev_{L/K}$ is asymptotically split if and only if
$\ev_{L/K}(i)$ asymptotically split for any $i\in I$. Similarly
$\ev_{L/K}$ is asymptotically split if and only if $\ev_{L/K}(i)$
is asymptotically split for any $i\in I$.
\end{prop}
\begin{proof}We will only give a proof for $\ev_{L/K}$; the dual
case is treated completely analogously. Assume that each of the
maps $\ev_{L/K}(i)$ is asymptotically split and let $W\in A\otimes
L$ be a finite dimensional $L$-submodule in $A$. We want to show
that the composite map
\[\xymatrix{W\ar[r]&(\bigotimes_iA_i\otimes
L)=A\otimes L\ar^-{\ev_{L/K}}[r]&\Map(\prod_iA_i(K),L})\] splits.
Without loss of generality we could assume that $W$ is of the form
$\bigotimes_i^nW_i\otimes L$ where $W_i$ is a finite dimensional
free $k$-submodule in $A_i$ for $i=1,2,\ldots, n$.

Note that \[\Map(A_i(K),L)\cong \lim_{\leftarrow}\Map(X_i,L)\cong
\] where the inverse limit is taken
over all finite subsets $X_i\in A_i(K)$ for $i=1,2$. Since there
exists a continuous splitting $\Map(A_i(K),L)\lra W_i$ and $W_i$
is finite dimensional it follows that this splitting map factors
through $\Map(X_i,L)$ for some finite subset $X_i\in
\Map(A_i(K),L)$. So we showed that there exist maps
$g_i:\Map(X_i,L)\lra W_i$ which split the compositions
\[W_i\lra A_i\otimes L\lra \Map(A_i(K),L)\lra \Map(X_i,L).\] Then
the composition \[\xymatrix{\Map(\prod_i A_i(K),L)\ar[r]&
\Map(\prod_i X_i,L)\cong
\bigotimes_i\Map(X_i,L)\ar[r]^-{\otimes_ig_i}&\bigotimes
_iW_i=W}\] is the desired splitting.

Conversely, assume that $\ev_{L/K}$ is asymptotically split and
let $W_i\in A_i\otimes L$ be an arbitrary finite dimensional
$L$-submodule. The the inclusion
\[W_i\hookrightarrow A_i\otimes L\rightarrow \bigotimes_iA_i\otimes L\lra
\Map(\prod_iA_i(K) ,L)\] splits and it follows that the inclusion
\[W_i\hookrightarrow A_i\otimes L\lra
\Map(A_i(K),L)\] also splits.
\end{proof}

\section{Numerical functions and functionals}
\begin{defi}
Let $\hat{L}$ be the field of fractions of $L$. The ring of
numerical functions $Num_{L/K}(A)$ is defined from the pullback
diagram of rings:
\[\xymatrix{Num_{L/K}(A)\ar[r]\ar[d]&\Map(A(K),L)\ar[d]\\
A\otimes\hat{L}\ar[r]^{\ev_{\hat{L}/K}}&\Map(A(K),\hat{L})}\] Here
the right vertical map is induced by the inclusion $L\lra
\hat{L}$. The $L$-module of numerical functionals is defined in a
dual fashion as a pullback of the following diagram:
\[\xymatrix{\Num_{L/K}(A)\ar[r]\ar[d]&\Hom(A,L)\ar[d]\\
\hat{L}[A(K)]\ar[r]^{\ev^*_{\hat{L}/K}}&\Hom(A,\hat{L})}\]
\end{defi}
\begin{rem}
If $A$ is a Hopf algebra over $K$ then $\Num_{L/K}(A)$ is an
$L$-algebra.
\end{rem}
It follows directly from the definition that the evaluation map
$\ev_{L/K}$ factors through numerical functions so that the
following diagram of algebras is commutative:
\begin{equation}\label{factor1}\xymatrix{A\otimes L\ar^{\ev_{L/K}}[rr]\ar[dr]^{e_1}&&\Map(A(K),L)\\
&Num_{L/K}(A)\ar[ur]^{e_2}}\end{equation} where the $e_1$ and
$e_2$ are natural inclusions. Dually we have the following
commutative diagram of $L$-modules:
\begin{equation}\label{factor2}\xymatrix{L[A(K)]\ar^{\ev^*_{L/K}}
[rr]\ar[dr]^{\underline{e}_2}&&\Hom(A,L)\\
&\Num_{L/K}(A)\ar[ur]^{\underline{e}_1}}\end{equation} Finally,
taking the (continuous) $L$-dual to the diagram (\ref{factor2})
and splicing it with (\ref{factor1}) we get the following two
factorizations of the map $\ev_{L/K}$:
\begin{equation}\label{factor3}\xymatrix{&\Num_{L/K}^*(A)\ar_{\underline{e}^*_2}[dr]\\
A\otimes
L\ar_{\underline{e}^*_1}[ur]\ar^{\ev_{L/K}}[rr]\ar[dr]^{e_1}&&\Map(A(K),L)\\
&Num_{L/K}(A)\ar[ur]^{e_2}}\end{equation} The similar
factorization of $\ev^*_{L/K}$ has the following form:
\begin{equation}\label{factor4}\xymatrix{&Num_{L/K}^*(A)\ar_{{e}^*_1}[dr]\\
L[A(K)]\ar_{{e}^*_2}[ur]\ar^{\ev_{L/K}}[rr]\ar[dr]^{\underline{e}_2}&&\Hom(A,L)\\
&\Num_{L/K}(A)\ar[ur]^{\underline{e}_1}}\end{equation}

\begin{example} Let $A=\mathbb{Z}[x]$. In this important case we
will shorten the notation $Num_{L/K}(A)$ to $Num_{L/K}$ and call
it the ring of \emph{$L/K$-numerical polynomials}. Furthermore
assume that $K=L=\mathbb{Z}$. In this case
$Num_{\mathbb{Z}/\mathbb{Z}}$ is the well-known ring of numerical
polynomials which consists of those polynomials with rational
coefficients which assume integer values at integers. It is
additively generated by the polynomials
\[P_k(x)={x\choose k}=\frac{x(x-1)\ldots(x-k+1)}{k!}\] for
$k=0,1.2,\ldots.$

Furthermore, the ring $Num_{\mathbb{Z}/\mathbb{Z}}$ is a Hopf
subalgebra inside $\mathbb{Q}[x]$. Its integral dual
$Num^*_{\mathbb{Z}/\mathbb{Z}}$ is isomorphic as a ring to the
ring of formal power series $\mathbb{Z}[[y]]$. For topologists we
offer the following explanation: recall that
$Num_{\mathbb{Z}/\mathbb{Z}}$ is isomorphic to
$KU_*(\mathbb{C}P^\infty)\otimes_{KU_*}\mathbb{Z}$, cf. \cite{Ba}.
Then by the universal coefficients formula
\[Num^*_{\mathbb{Z}/\mathbb{Z}}\cong
\Hom_{KU_*}(KU_*(\mathbb{C}P^\infty),\mathbb{Z})\cong
KU^*(\mathbb{C}P^\infty)\otimes_{KU_*}\mathbb{Z}\cong
\mathbb{Z}[[y]].\] The upper part of the diagram (\ref{factor4})
takes the following form:
\[\xymatrix{&\mathbb{Z}[[y]]\ar_{{e}^*_1}[dr]\cong Num_{\mathbb{Z}/\mathbb{Z}}^*\\
\mathbb{Z}[t,t^{-1}]\cong\mathbb{Z}[\mathbb{Z}]\ar_{{e}^*_2}[ur]\ar^{\ev_{L/K}}[rr]&&\mathbb{Z}[x]^*
\cong\Gamma_{\mathbb{Z}}}\] Direct inspection shows that
$e_2^*(t)=y+1$ and $e_1^*(y)=\sum_{i=1}^\infty u_i\in
\Gamma_\mathbb{Z}$.

\end{example}
The following result is
a useful criterion for (non)existence of nontrivial numerical
polynomials.   Related results could be found in \cite{CC}.
\begin{prop}\label{invert}
If there exists an element $a\in L$ such that $a$ is not
invertible in $L$ and for which the image of the composite map
\[\pi_a:K\lra L\lra L/a\] is finite then $Num_{L/K}\neq L[x]$. Conversely, suppose
$Num_{L/K}\neq L[x]$ and, in addition,  $L$ is a unique
factorization domain (UFD). Then there exists a noninvertible
element $a\in L$ for which $\pi_a$ has a finite image.\end{prop}
\begin{proof}
Suppose that there exists $a\in L$ for which $\im(\pi_a)$ consists
of $m$ elements $e_1,\ldots, e_m$. Choose a collection of
representatives $k_i\in K$ of these residue classes and consider
the polynomial
\[P(x)=\frac{\prod_{i=1}^m(x-k_i)}{a}\in \hat{L}[x].\]
Note that since $a$ is not invertible in $L$ the polynomial $P(x)$
is not contained in $L[x]$. Let $k\in K$. Then the collection
$\{k-k_i\}$ contains an element congruent to zero module $a$.
Therefore $P(k)\in L$.

Conversely, suppose that $P(x)\in Num_{L/K}$ and $P(x)$ is not
contained in $L[x]$. Then $P(x)=\frac{Q(x)}{a}$ with $Q(x)\in
L[x]$, the greatest common divisor of all coefficients of $Q(x)$
is a unit in $L$ and $a\in L$ is noninvertible in $L$. Without
loss of generality we assume that $a$ is prime in $L$. We claim
that $\im \pi_a$ is finite.  Indeed, if this is not the case there
exists an infinite collection $\{k_i,i=1,2,\ldots\}\in K$ of
representatives of $\im\pi_a$ in $K$ which are pairwise distinct.

Denote by $\bar{Q}(x)$ the image of $Q(x)$ under the map $L[x]\lra
L/a[x]$. Then $\bar{Q}[x]$ is not identically zero and
$\bar{Q}(k_i)=0$ for all $i$. In other words $\bar{Q}(x)$ has
infinitely many roots which is impossible since $L/a$ is an
integral domain (recall that $a$ is a prime element in $L$). This
contradiction finishes the proof.
\end{proof}
Rings of integers in algebraic number fields and algebraic
function fields provide natural examples of nontrivial (i.e.
having nontrivial denominators) numerical polynomials. We restrict
ourselves with giving two examples.
\begin{example}
$K=L=\f_p[q]$. Then the polynomials
\[f_n(x)=\frac{x^{p^n}-x}{q^{p^n}-q}\] belong to $Num_{L/K}$.
\end{example}
\begin{example}
$K=L=\mathbb{Z}[i]$, the ring of Gaussian integers. A nontrivial
example of $L/K$-numerical polynomials is given by
\[f_n(x)=\frac{1}{n!}\prod_{0\leq a,b<n}(x-a-ib).\]
\end{example}

Now we will consider the notion of numerical functionals in more
detail. Even in the case $K=L=\mathbb{Z}, A=\mathbb{Z}[x]$ they
have not been studied to the best of our knowledge. To give
nontrivial examples of $\Num_{L/K}(A)$ suppose, as before that
$A=K[x]$ in which case we will use the notation $\Num_{L/K}$.
Recall that $\Num_{K/K}$ is a subring in the group ring
$\hat{L}[K]$.

Let $\omega:=(k_1,\ldots,k_n)$ be a collection of pairwise
distinct elements in $K$. Let $A_\omega$ be the classical
Vandermonde matrix:
\begin{equation}\label{van}A_{\omega}= \begin{pmatrix}1&1&\ldots&1\\
k_1&k_2&\ldots&k_{n}\\
\cdot&\cdot&\ldots &\cdot\\k_1^{n-1}&k_2^{n-1}&\ldots
&k_{n}^{n-1}\end{pmatrix}\end{equation}  Set $d_\omega:=\det
A_\omega$. Further denote by $e_{l,\omega}$ the element in
$\hat{L}[K]$ defined as $\frac{1}{d_{\omega}}\det \tilde{A}$ where
$\tilde{A}$ is the matrix obtained from the matrix $A_{\omega}$ by
replacing formally the $l$th row by $([k_1],[k_2],\ldots,[k_n])$.
\begin{example}
Let $\omega=(k_1,k_2)$. Then $d_{\omega}=k_2-k_1$ and
\[e_{1,\omega}=\frac{1}{k_2-k_1}\det\begin{pmatrix} [k_1]&[k_2]\\k_1&k_2\end{pmatrix}
=\left(\frac{k_2}{k_2-k_1}\right)[k_1]+
\left(\frac{k_1}{k_1-k_2}\right)[k_2];\]
\[e_{2,\omega}=\frac{1}{k_2-k_1}\det\begin{pmatrix} 1&1\\ [k_1]&[k_2]\end{pmatrix}
=\left(\frac{1}{k_1-k_2}\right)[k_1]+
\left(\frac{1}{k_2-k_1}\right)[k_2].\]
\end{example}
\begin{theorem}\label{collection}\begin{enumerate}\item
The elements $e_{l,\omega}$ belong to $\Num_{L/K}$ for any
$\omega,l$. Moreover,   $\Num_{L/K}$ is generated by
$e_{l,\omega}$ as an $L$-module.\item The canonical inclusion
$\underline{e}_1^*:\Num_{L/K}\lra\Gamma_L$ has a dense
image.\end{enumerate}
\end{theorem}\begin{proof}
Denote by $A_{l,m,\omega}$ the matrix obtained from $A_\omega$ by
replacing the $l$th row with $(k_1^m,k_2^m,\ldots,k_n^m)$. Note
that $A_{l,l-1,\omega}=A_\omega$ for all $l=1,\ldots,n$. Set
$d_{l,m,\omega}:=\det A_{l,m,\omega}$.

The straightforward computation using elementary properties of
determinants yields the following key formula.
\[
\ev_{L/K}^*(e_{l,\omega})=\sum_{m=0}^\infty\frac{d_{l,m,\omega}}{d_\omega}u_m.\]
Since $d_{l,m,\omega}=0$ for $0 \leq m \leq n-1, m \neq l$ and
$d_{l,l-1,\omega}=d_\omega$ we see that
\begin{equation}\label{important}\ev_{L/K}^*(e_{l,\omega})=u_{l-1}+\sum_{m=n}^\infty
\frac{d_{l,m,\omega}}{d_\omega}u_m,\end{equation}

The expression $\frac{d_{l,m,\omega}}{d_\omega}, m \geq n,$ is (up
to sign) the classical Schur function in $k_1,\ldots,k_n$ (cf.
\cite{macdonald}) corresponding to the partition
$(m-n+1,0,\ldots,0)$ of length $n$. The Schur functions are
symmetric polynomials with integer coefficients. That shows that
$e_{l,\omega}\in \Num_{L/K}$ for all $l,\omega$.

Furthermore, formula (\ref{important}) also implies that the image
of $e_1^*$ is dense in $\Gamma_L$. We still need to prove that
$e_{l,\omega}$ generate $\Num_{L/K}$. Let
$\Gamma_L^n:=\Gamma_L/(u_i,i\geq n)$, the quotient of $\Gamma_L$
by the submodule spanned by $u_n, u_{n+1},\ldots$. Set
$V_\omega:=\hat{L}\langle[k_1],\ldots,
[k_n]\rangle\bigcap\Num_{L/K}$. Composing $\ev^*_{\hat{L}/K}$ with
the projection $\Gamma_L\lra\Gamma_L^n$ we get the map
$\phi^n_\omega:V_\omega\lra\Gamma^n_L$. Let us prove the following
lemma. \begin{lem}\label{include} The map $\phi^n_\omega$ is a
monomorphism. \end{lem} \begin{proof} Let
$f=\sum_{i=1}^n\xi_i[k_i]$ where $\xi_i\in\hat{L}$. Then the
condition $\phi_\omega^n(f)=0$ implies that
$\sum_{i=1}^n\xi_ik_i^m=0$ for all $m=0,1,\ldots,n-1$. Since the
relevant Vandermonde matrix is nondegenerate we conclude that
$\xi_i=0,i=1,2,\ldots,n$. \end{proof} Consider the $L$-module
$L\langle e_{1,\omega},e_{2,\omega},\ldots,e_{n,\omega}\rangle$,
the $L$-submodule in $V_\omega$ spanned by
$e_{1,\omega},e_{2,\omega},\ldots,e_{n,\omega}$. From formula
(\ref{important}) we conclude that $\phi_\omega^n$ maps $L\langle
e_{1,\omega},e_{2,\omega},\ldots,e_{n,\omega}\rangle$
isomorphically onto $\Gamma^n_L$. Combining this with Lemma
\ref{include} we see that $\phi^n_\omega:V_\omega\lra\Gamma^n_L$
is an isomorphism and therefore $L\langle
e_{1,\omega},e_{2,\omega},\ldots,e_{n,\omega}\rangle=V_\omega$.
Thus, the elements $e_{l,\omega}$ span the whole $\Num_{L/K}$ as
claimed. \end{proof} \begin{rem} Note that our proof shows that
$\Num_{L/K}$ is a union of the submodules $L\langle
e_{1,\omega},e_{2,\omega},\ldots,e_{n,\omega}\rangle$ and each of
these submodules maps isomorphically onto
$\Gamma_L^n=\Gamma_L/(u_i,i\geq n)$ under the composition
\[\Num_{L/K}\lra \Gamma_L\lra\Gamma^n_L.\]\end{rem}
\begin{cor}\label{equiv} The ring of numerical functionals
$\Num_{L/K}$ coincides with $L[K]$ if and only if $L$ contains
$\hat{K}$, the field of fractions of $K$. \end{cor} \begin{proof}
By Theorem \ref{collection} the element $\frac{[k]-[0]}{k}\in
\Num_{L/K}$ for any $k\in K$. Assume that $\Num_{L/K}=L[K]$. Since
$[k]-[0]$ could be taken to be one of the basis elements in the
free $L$-module $L[K]$ we conclude that $k$ is invertible in $L$.

Conversely, suppose that $L$ contains $\hat{K}$. From the
construction of the basis elements $e_{l,\omega}\in \Num_{L/K}$ it
follows that $e_{l,\omega}\in L[K]$. Therefore $\Num_{L/K}=L[K]$
in this case. \end{proof} \begin{rem}As an aside mention that
numerical functionals could be considered as difference operators
with constant coefficients. Namely, suppose for simplicity that
$K=L=\mathbb{Z}$. We associate to an element
$\omega=\Sigma_i\alpha_i[k_i]\in \mathbb{Q}[\mathbb{Z}]$ the
operator $\bar{\omega}$ acting on $\mathbb{Q}[x]$ according to the
rule \[p(x)\mapsto \Sigma_i\alpha_ip(x+k_i)\] where $p(x)\in
\mathbb{Q}[x]$. Then clearly if $\bar{\omega}$ maps
$\mathbb{Z}[x]\in \mathbb{Q}[x]$ into itself then $\omega$ is a
numerical functional. Conversely, if $\omega$ is a numerical
functional then it is easy to see, using Taylor's expansion
formula that $\bar{\omega}$ restricts to an operator on
$\mathbb{Z}[x]$. \end{rem}
 \section{Properties of the
evaluation map and its dual for polynomial rings} We will now
consider the properties of the evaluation map and its dual in the
case $A=K[x]$ from the point of view of numerical functions and
numerical functionals. Informally speaking, the existence of
nontrivial numerical polynomials (i.e. those not having
coefficients in $L$) is an obstruction to the map
$\ev^*_{L/K}:L[x]\lra \Map(K,L)$ being an asymptotically split
monomorphism. Dually, the existence of nontrivial numerical
functionals (i.e. polynomials in $\Num_{L/K}$ but not in
$\Gamma_L$ implies that $\ev_{L/K}:L[K]\lra \Gamma_L$ is not
asymptotically split.

More precisely, we have the following result.
\begin{theorem}\label{main1} The following conditions are
equivalent. \begin{enumerate}\item There exist no nontrivial
numerical functionals, i.e. $\Num_{L/K}=L[K]$.\item The map
$\ev^*_{L/K}:L[K]\lra \Gamma_L$ is asymptotically split
(equivalently the map $\ev_{L/K}:L[x]\lra \Map(K,L)$ has a dense
image). \item The ring $L$ contains $\hat{K}$, the ring of
fractions of $K$. \end{enumerate} \end{theorem} \begin{proof}
First observe that the equivalence of (1) and (3) is the statement
of Corollary \ref{equiv}. Further assume that
$\ev^*_{L/K}:L[K]\lra \Gamma_L$ is asymptotically split and
consider $f\in \Num_{L/K}\subset \hat{L}[K]$. There exists $l\in
L$ for which $lf\in L[K]$. It follows that the inclusion $L\langle
lf\rangle\hookrightarrow Num_{L/K}$ splits where $L\langle
lf\rangle$ is a one-dimensional $L$-submodule in $L[K]$ spanned by
$lf$. Therefore $f\in L[K]$. Thus we proved (2)$\Rightarrow$(1).

Let us prove (3)$\Rightarrow$(2). Suppose that $L$ contains
$\hat{K}$ and let $W$ be a finite dimensional free $L$-submodule
in $L[K]$. Without loss of generality we assume that $W$ is
spanned by the finite collection of elements $[k_i]\in K$ and
denote by $W_{\hat{K}}$ the $\hat{K}$-submodule in $\hat{K}[K]$
spanned over $\hat{K}$ by $[k_i]$. Recall that $\Gamma^i_L$ is the
quotient $\Gamma_L/(u_k,k\leq i)$ and similarly denote
$\Gamma^i_{\hat{K}}:=\Gamma_{\hat{K}}/(u_k,k\leq i\rangle)$. Since
$W_{\hat{K}}$ is finite dimensional over $\hat{K}$ there exists a
positive integer $i$ for which the composition
\[f:W_{\hat{K}}\hookrightarrow \hat{K}[K]\lra\Gamma_{\hat{K}} \lra
\Gamma^i_{\hat{K}}\] is an inclusion. Then there exists a
splitting map $g:\Gamma^i_{\hat{K}}\lra W_{\hat{K}}$ such that
$g\circ f=id_{W_{\hat{K}}}$. Since $\Gamma^i_{\hat{K}}$ is finite
dimensional over $\hat{K}$ we have $\Gamma^i_{L}\cong
\Gamma^i_{\hat{K}}\otimes_{\hat{K}}L$. Therefore tensoring the
splitting map $g$ with $L$ over $\hat{K}$ we obtain a map
$\hat{g}:\Gamma^i_L\lra W$ which splits the composite map
\[W\hookrightarrow L[K]\lra\Gamma_L \lra \Gamma^i_L.\]
Precomposing $\hat{g}$ with the projection $\Gamma_L \lra
\Gamma^i_L$ we obtain a continuous map $\Gamma_L\lra W$ which
splits the inclusion $W_\hookrightarrow L[K]\lra\Gamma_L$ as
desired. \end{proof} \begin{example}Let $K=\mathbb{Z}$ and
$L=\mathbb{Q}$. The map
\[ev_{\mathbb{Q}/\mathbb{Z}}^*:\mathbb{Q}[\mathbb{Z}]\lra\Gamma_{\mathbb{Q}}\]
is an asymptotically split monomorphism.\end{example} We will now
consider the question when the map $\ev_{L/K}$ is asymptotically
split. In preparation for this let us introduce the notion of the
ring of \emph{DS-extension}. Here DS is an abbreviation for
`divisible sequence'. \begin{defi}Let $Q$ be a subring of a
commutative ring $R$ with no zero divisors. We say that $R$ is a
DS-extension of $Q$ if there exists an infinite sequence of
elements $q_i\in Q, i=1,2,\ldots$ such that all differences
$q_i-q_j$ are invertible in $R$ for $i\neq j$. \end{defi} We will
now give various examples of DS-extensions. \begin{enumerate}
\item If $Q$ is an infinite ring and $R$ contains the field of
fractions of $Q$, then $R$ is a DS-extension of $Q$. \item A local
ring with infinite residue field is a DS-extension of itself.
\item The so-called Novikov ring $\mathbb{Z}((q))$ consisting of
Laurent power series with integer coefficients is a DS-extension
of its subring $\mathbb{Z}[q]$\item A more economical variant of
the previous example is $\mathbb{Z}[x^{\pm 1}][\{(x^n-1)^{-1}\}]$,
$n=1,2,\ldots$, the ring of polynomials where the elements $x$ and
all cyclotomic polynomials are inverted.\item Any DS-extension of
$\mathbb{Z}$ must contain $\mathbb{Q}$. Indeed, suppose that
$k_i,i=1,2,\ldots$ are integers such that $k_i-k_j$ is invertible
in some ring $R\supset\mathbb{Z}$ for $i\neq j$. Since the
sequence $\{k_i\}$ is infinite for any prime number $p$ there
exists a pair $i,j$ with $i\neq j$ and $k_i\equiv k_j\mod p$.
Therefore $p$ must be invertible in $R$ and $\mathbb{Q}\subset R$.
Similarly if $Q$ is a principal ideal domain such that for any
maximal ideal $(q)\in Q$ the residue field $Q/(q)$ is finite then
any DS-extension of $Q$ must contain $\hat{Q}$, the field of
fractions of $Q$. \end{enumerate} We can now formulate our second
main theorem which is a partial dualization of Theorem \ref{main1}
\begin{theorem}\label{main2} We have the following chain of
implications for the conditions $(1), (2), (3), (4)$ below:
$1)\Rightarrow (2)\Rightarrow (3)\Rightarrow (4)$. If $L$ is a UFD
then $(4)$ implies $(3)$. If, further, $L$ is a discrete valuation
ring (DVR) then $(4)$ implies $(1)$. \begin{enumerate}\item $L$ is
a $DS$-extension of the ring $K$.\item The map $\ev_{L/K}:L[x]\lra
\Map(K,L)$ is asymptotically split (equivalently map
$\ev^*_{L/K}:L[K]\lra \Gamma_L$ has a dense image)\item There
exist no nontrivial numerical polynomials, i.e.
$Num_{L/K}=L[x]$.\item For any element $a\in L$ the composite map
$K\hookrightarrow L\lra L/(a)$ has an infinite image.
\end{enumerate} \end{theorem} \begin{proof}(1)$\Rightarrow$(2).
 Let $L$ be a DS-extension of $K$
Recall that $\Gamma_L$ has a topological basis $u_n$.  It is clear
that the image of $\ev^*_{L/K}$ is dense in $\Gamma_L$ if in
$L[K]$ there exists a collection of elements
$a_1,\ldots,a_n,\ldots$ such that
\[\ev^*_{L/K}(a_n)=u_n+\sum_{i\geq 1}\beta_{n,i}u_{n+i}.\] Set
$a_n=\sum_{j=1}^n\alpha_{n,j}[k_j]$.
Then
\[ev^*_{L/K}(a_n)=\sum_{l=0}^\infty(\sum_{j=1}^n\alpha_{n,j}k^l_j)u_l.\]
Therefore, one needs to find a collection of elements
$\alpha_{n,j}\in L$ and $k_j\in K$ where $j=1,\ldots, n$ such that
the following equations hold for $l=0,\ldots,n-1$.
\begin{align}\label{vande}\nonumber\sum_{j=1}^n\alpha_{n,j}k^l_j=0;\\
 \sum_{j=1}^n\alpha_{n,j}k^n_j=1.\end{align}
Since $L$ is a $DS$-extension of $K$ there exists a system of
elements $k_j$ where $j=0,\ldots,n$ and such that all differences
$k_i-k_j$ are invertible in $L$ for $i\neq j$. Then the matrix of
(\ref{vande}) considered as a linear system with respect to the
unknowns $\alpha_{n,j}$ is the well-known Vandermonde matrix and
its determinant equals $\prod_{i<j}(k_i-k_j)\neq 0$. Therefore
this linear system admits a solution.

(2)$\Rightarrow$(3). This implication is almost obvious and proved
along the same lines as (2)$\Rightarrow$(1) in the proof of
Theorem \ref{main1}.

(3)$\Rightarrow$(4). This is proved in Proposition \ref{invert}.
The same proposition provides a partial inverse
(4)$\Rightarrow$(3) provided that $L$ is a UFD.

Finally let us prove the implication (4)$\Rightarrow$(1) assuming
that $L$ is a DVR. Let $\pi\in L$ be the generator of the maximal
ideal of $L$. Choose the infinite sequence $k_i,i=1,2,\ldots$ of
distinct elements in the image of the composite map
$K\hookrightarrow L\lra L/(\pi)$. Let $\tilde{k}_i\in L$ be an
arbitrary representative of residue class $k_i\in L(\pi)$. Then
for any $i\neq j$ the difference $\tilde{k}_i-\tilde{k}_j$ does
not belong to the ideal $(\pi)$ and therefore is invertible. Thus,
$L$ is a DS-extension of $K$.
\end{proof}
\begin{rem}\label{several}
Theorems \ref{main1} and \ref{main2} give criteria for the
evaluation map and its dual to be asymptotically split (or to have
dense images) in the case of the polynomial ring in one variable.
Using Proposition \ref{tensor} we see that the same criteria can
be used for the rings of polynomials in several variables.
\end{rem}

\section{Applications to algebras of cohomology operations}
In this section we show that the Landweber-Novikov and Steenrod
algebras after a certain extension of scalars and completion
become isomorphic to the (completions of) some group rings. This
is analogous to the Theorem 6.2.3 of \cite{Ra} saying that the
continuous dual to the Morava stabilizer algebra over
$\mathbb{F}_{p^n}$ is isomorphic to a group algebra. Related
results in the context of the classical Steenrod algebra are
contained in the recent preprint of J. Palmieri, \cite{Pal}.

Let $R$ be a commutative ring. The set of formal power series of
the form
\[x(t)=t+x_1t^2+\ldots+x_kt^{k+1}+\ldots\]
forms a group under composition of formal power series. This group
plays an important role in symplectic geometry, algebraic
topology, singularity theory and other fields. In group theory it
is known by the name \emph{Nottingham group}  cf. \cite{cam}. It
will be denoted by $N(R)$.

 The group $N(R)$ is the group of
$R$-points of a certain Hopf algebra $N$. This Hopf algebra is
isomorphic to the ring $\mathbb{Z}[x_1,x_2,\ldots]$ as an algebra.
The diagonal $\Delta$ is dual to the composition of the power
series and is given by the following formula:
\[\Delta(x_k)=\sum_{i+j=k}x_i\otimes(t+x_1t^2+x_2t^3+\ldots)^{j+1}_i\]
where $(t+x_1t^2+x_2t^3+\ldots)^{j+1}_i$ denotes the coefficient
at the $i$th power of $t$ in the series
$(t+x_1t^2+x_2t^3+\ldots)^{j+1}_i$.

An important property of the Hopf algebra $N$ is that it is
graded. Indeed, letting $|x_i|=i$ we see that the diagonal in $N$
preserves grading. Moreover, each graded component of $N$ a is a
free abelian group of finite rank.

The Hopf algebra $N$ is closely related to \emph{complex cobordism
theory} $MU^*$. Namely, the algebra of formal differential
operators on $N$ is isomorphic to the algebra of cohomology
operations in the $MU^*$. Under this isomorphism the subalgebra of
left-invariant differential operators corresponds to the so-called
\emph{Landweber-Novikov} algebra $S$. We refer the reader to the
paper by the first author \cite{Buch} for an exposition of these
results. We would like to take another approach, not using
differential operators.

\begin{defi}The Landweber-Novikov algebra $S$ is the graded
$\mathbb{Z}$-dual of the Hopf algebra $N$. The complete
Landweber-Novikov algebra $\hat{S}$ is the \emph{ungraded}
$\mathbb{Z}$-dual of the Hopf algebra $N$.
\end{defi}

Observe that $\hat{S}$ is a topological Hopf algebra over
$\mathbb{Z}$. For a ring $R$ the algebra $\hat{S}\otimes R$ has a
linear topology inherited from $\hat{S}$ and we denote by
$S\hat{\otimes}R$ the completion of $\hat{S}\otimes R$ with
respect to this topology. Note that $\hat{S}\hat{\otimes}R\cong
\Hom(N,R)$.

We have the following result.
\begin{prop}\label{landweber}Let $R$ be a $DS$-extension of a ring $Q$.
Then the topological Hopf $R$-algebra $\hat{S}\hat{\otimes} R$,
has a dense Hopf $R$-subalgebra isomorphic to $R[N(Q)]$.\end{prop}

 There is an analogue of this result in the context
of the Steenrod algebra. Namely, consider for any $\f_p$-algebra
$R$ the group $P(R)$ consisting of all formal power series of the
form $\sum_{i=0}^\infty a_ix^{p^i}$ where $a_i\in R$ and $a_0=1$.
The group operation is the composition of power series.
(Incidentally, $P(R)$ is the group of strict automorphisms of the
additive formal group law over $R$.) Then the functor $R\mapsto
P(R)$ is represented by a Hopf algebra $P=\f_p[x_1,x_2,\ldots]$
where the diagonal map is given by the formula \[\Delta x_n
=x_n\otimes 1+1\otimes x_n+\sum_{k=1}^\infty x_k\otimes
y_{i-k}^{p^k}.\] The Hopf algebra $P$ is in fact graded with
$deg(x_i)=p^i-1$. The graded $\f_p$-dual to $P$ will be denoted by
$\mathcal{A}$. The algebra $\mathcal{A}$ is isomorphic to the
Steenrod algebra if $p=2$ and to the algebra of reduced powers if
$p$ is odd. We will use the symbol $\hat{\mathcal{A}}$ to denote
the ungraded dual to $P$; it is a topological Hopf algebra over
$\f_p$. For any ring $F$ denote by
$\hat{\mathcal{A}}\hat{\otimes}F$ the completed tensor product of
$\hat{\mathcal{A}}$ and $F$, clearly
$\hat{\mathcal{A}}\hat{\otimes}F\cong \Hom(P,F)$.

Then the analogue of Proposition \ref{landweber} reads as follows:
\begin{prop}\label{steenrod}Let $R$ be a $DS$-extension of a ring $Q$ and suppose that $Q$
is an $\f_p$-algebra. Then the topological Hopf $R$-algebra
$\hat{\mathcal{A}}\hat{\otimes} R$, has a dense Hopf
$R$-subalgebra isomorphic to $R[P(Q)]$
\end{prop}
The statements of propositions \ref{landweber} and \ref{steenrod}
follow from Theorem \ref{main2}, taking into account Remark
\ref{several}.

\end{document}